\documentclass[11pt]{article} 
\usepackage[utf8]{inputenc} 
\usepackage[polutonikogreek,english]{babel}
\usepackage[or]{teubner}
\usepackage{xpatch}
\xapptocmd{\greektext}{\edef~{\string~}}{}{}
\usepackage{amsfonts}
\usepackage{amsmath}
\usepackage{amsthm}
\usepackage{titlesec}

\titleformat{\section}  
  {\fontsize{13}{15}\bfseries} 
  {\thesection.} 
  {1em} 
  {} 
  [] 

\titleformat{\subsection}  
  {\fontsize{12}{14}\bfseries} 
  {\thesubsection.} 
  {1em} 
  {} 
  [] 

\titleformat{\subsubsection}  
  {\fontsize{11}{13}\bfseries} 
  {\thesubsubsection.} 
  {1em} 
  {} 
  [] 

       \usepackage[polutonikogreek,english]{babel}

\usepackage{hyperref}
\usepackage{geometry} 
\usepackage{graphicx}
\graphicspath{ {c:/users/fpc-dprotopapas/documents/latex/} }

\usepackage{booktabs} 
\usepackage{array}
\usepackage{paralist} 
\usepackage{verbatim} 
\usepackage{subfig}
\usepackage{multirow}

\usepackage{fancyhdr} 
\usepackage[nottoc,notlof,notlot]{tocbibind} 
\usepackage[titles,subfigure]{tocloft}

\usepackage{epigraph}
\title{Bernhard Riemann: Reflections on a Bicentennial}
 \author{Athanase Papadopoulos}

\date{\today}

\begin{document}

     \maketitle
   
%
%
%
%
%
%

\epigraph{\itshape 
 After a quiet preparation Riemann came forward like a bright meteor, only to be extinguished
soon afterwards.

(Felix Klein, Development
of mathematics in the 19th century, 1926, English transl. R. Hermann,)}
   
   \begin{abstract}

   This paper was written at the occasion of the 
200th anniversary of the birth of Bernhard Riemann (1826-1866). Within a limited framework, we have tried to include, in addition to a few chosen elements of Riemann's mathematical work, aspects of his ideas on physics and philosophy. The focus is also on Riemann’s character, the direct influences he was subject to, and the impact of his work on some chosen group mathematicians that followed him. I have included this in the broader idea of historical continuity in mathematics. 
The presentation is based on Riemann’s correspondence with his family and other firsthand accounts about his life and work. 
   The final version of this article will appear in the Indian Mathematics Consortium
Bulletin.  
   The paper was also the subject of a colloquium talk I gave at the Chebyshev laboratory of the Saint Petersburg State University, in October 2026.
   
   \medskip
   
   \medskip

   
   \end{abstract}

\section{Introduction: Historical continuity in mathematics}
This year, we celebrate the
200th anniversary of the birth of Bernhard Riemann (1826-1866). I would like to start my commemoration by saying a few general words about our mathematical ancestors, who bequeathed to us the theories we are building upon, the problems we are working on, and the tools to solve them.

More so than in the other sciences, historical continuity is  particularly evident in mathematics. Every mathematician contributes, to the best of his ability, to explore this coherent and hermetic universe of mathematics --- hermetic because only those who work within it can contemplate it and feel its secrets. In mathematics, and unlike in the other sciences, fundamental concepts and  theories never die; they are passed down from generation to generation, and it often happens that each receiving generation can use them with almost no change. One need only think of Euclid’s \emph{Elements}, where the axioms, definitions, and propositions are still in use  about 2300 years after they were written, with the very same words. Euclid is not an isolated case: Nearly two-thousand-year-old treatises by Archimedes, Apollonius, Diophantus, Theodosius, Menelaus and others could serve as bases for courses taught to our students today --- and sometimes they do; the author of these lines gave courses on non-Euclidean geometry based on Menelaus's \emph{Spherics} \cite{RP-Menelaus}, a work written in the first century of our era. The editor of an English modern version of Leonhard Euler's \emph{Introductio in analysin infinitorum} (Introduction to the analysis of the infinite)  writes that André Weil ``was trying to convince the mathematical community that our students of mathematics would profit much more from a study of Euler's \emph{Introductio in Analysin Infinitorum}, rather than of the available modern textbooks"  (see the Introduction to Vol. 1 of the English translation  \cite{Euler-Analysis}).  Mathematics, especially that written by the masters, never goes out of date.

Thus, it is not surprising that Riemann, in his writings, very frequently expresses his debt to his predecessors. For every topic on which he worked, he was thoroughly familiar with the existing literature. His 1854 Habilitation concise\footnote{The original publication is 18 pages long, with generous margins.} address, \emph{\"Uber die Hypothesen, welche der Geometrie zu Grunde liegen} (On the hypotheses which lie at the foundations of geometry), is full of references to his predecessors: Euclid, Legendre, Abel, Lagrange, Pfaff and others.
 In the final section of this address, we can read:
\begin{quote}\small
The progress of recent centuries in the knowledge of mechanics depends almost entirely on the exactness of the construction which
has become possible through the invention of infinitesimal calculus, and
through the simple principles discovered by Archimedes, Galileo, and Newton, and used by modern physics.\footnote{Regarding Riemann's  Habilitation, I am using the English translation by the mathematician and philosopher W. K. Clifford (1845-1879), see \cite{Riemann-Habilitation-Lecture}.}
\end{quote}
The references are eloquent: 
Archimedes, who developed the first integral calculus, allowing him to compute length, area and volume, Galileo, who introduced the
modern notions of motion, velocity and acceleration, and Newton, who
 was the first 
to give a mathematical expression to the forces of nature, which allowed him to describe
the motion of bodies in resisting media, and to whom is attributed a
celebrated notion of space, the ``Newtonian space".  
 
 As a matter of fact, the subject
of Riemann’s Habilitation lecture includes the three domains of Newton’s major \emph{opus}, the  \emph{Philosophi\ae  \ Naturalis Principia Mathematica (Mathematical Principles of Natural Philosophy)}:
mathematics, physics and philosophy. It is also interesting to note that Archimedes,
Galileo and Newton are quoted as the three founders of mechanics in the \emph{Discours préliminaire} (which plays the role of an Introduction) of Joseph Fourier’s \emph{Théorie analytique de la chaleur}  (Analytic theory of heat) \cite[p. i–ii]{Fourier-Heat}, a work in which the author lays down the foundations of the theory
of trigonometric series. Fourier's dissertation was written 32 years before Riemann's Habilitation thesis carrying the title \emph{\"Uber die Darstellbarkeit einer Function durch eine trigonometrische Reihe} (On the representability of a function by a trigonometric series) \cite{Riemann-Habilitation-memoir}, which is precisely
concerned with trigonometric series. Again, in this work, the care with which Riemann analyses
the evolution of the field he is expanding, since its birth, and the wealth of historical details he gives, are another indication of  how much he valued the history of mathematical ideas, and how much he was
aware of the historical developments of the topics he was working on. Like in his other writings, Riemann reviews in this dissertation the important paths and sometimes warns
against the wrong tracks that some of his predecessors took for the solutions of the various problems he was tackling. In the same work, he lays down the foundations of the modern theory of integration, after reviewing results and conjectures by Euler, d’Alembert, Lagrange,
Daniel Bernoulli, Dirichlet, Fourier and others on trigonometric series, a topic which served as a motivation for his development of what became known later as the Riemann integral.    It is also worth noting, since Riemann quotes Archimedes on several occasions, that Bourbaki considers the latter to be one of the main founders of integration theory. Let me quote Bourbaki's \emph{Elements of History of mathematics} \cite[p. 247-248]{Bourbaki-Elements-histoire}: \begin{quote}\small
Concerning integration, Cauchy’s work represents
a return to the sound traditions of antiquity and the early
17th century, but one based on technical methods that were still inadequate.
To define the definite integral, Cauchy returns to the method of exhaustion, or as we
would say, to ``Riemann sums" (which would be better called
``Archimedes sums", or ``Eudoxus sums").
\end{quote}
 One may recall here that Archimedes had developed the tools necessary for calculating areas of figures such as discs, sections of parabolas and sections of a sphere, considering analogues of what we call Riemann sums. On this matter, I would like to refer to a recent article by Demetrios Christodoulou \cite{Christodoulou}. As for Eudoxus, we owe to him the so-called \emph{exhaustion method}, used by Euclid in Book XII of the \emph{Elements}, a book which is due to Eudoxus. In fact, integration theory started precisely in Book XII. The exhaustion method continued developing and was greatly extended by Archimedes.

Riemann’s sense of history and historical continuity  is also plainly  manifest in the two-page announcement of his memoir
\emph{Beitr\"age zur Theorie der durch die Gauss'sche Reihe $F(\alpha, \beta, \gamma, x)$ darstellbaren
Functionen} (Contribution to the theory of functions representable by Gauss’s series
$F(\alpha, \beta, \gamma, x)$) published in 1857 (see \cite{Riemann-Anzeige-Beitrage} for the announcement and \cite{Riemann-Beitrage} for the memoir itself). In this memoir, he presents the origin of the questions on which he was working, related to the properties and behavior of hypergeometric functions, and while doing so, he mentions the names of Wallis, Euler, 
Pfaff, Gauss, Gudermann, Jacobi and Kummer. There are many other examples of Riemann paying tribute to his predecessors for their work on which he built.

 The next section contains, as every commemorative article should, some biographical details about Riemann. All along this section, I shall stress on his  teachers, another way to honor his predecessors.

  \section{\emph{Vita} and studies}

  As sources for Riemann's biography, we have at our disposal first-hand accounts that include letters he wrote to his family,
which allow us, to some extent, to reliably trace part of his life story. The letters are easily available, see e.g. \cite{Neuenschwander}, and \cite{Riemann-Letters} for English  translations of some of them.\footnote{Some dates differ from source to source and there may be mistakes in the transcription of the dates. We have done our best to transmit the content of some of these letters and to follow the natural chronology.}
Based on such documents, we can also get a sense of Riemann's feelings about his school studies, information on the courses he followed at university, his anxieties, and some of his needs. Together with notes about him written by his teachers, these documents give us some insight into his state of mind during the important period of his life when he was a school pupil, and then a university student. These letters and notes have helped various authors to write a number of biographies of Riemann.

 Georg Friedrich Bernhard Riemann, usually called simply Bernhard Riemann, was born on September 17, 1826, in Breselenz, a small isolated 
 farming and pastoral village of very few inhabitants, which was
   part of the Kingdom of Hanover, now in Germany. The region was poor and undeveloped. 
   His father, Friedrich Bernhard Riemann, was a Lutheran pastor of modest means. Soon after young Bernhard's  birth, the family moved to Quickborn, an even smaller village, situated about 15 kilometers from Breselenz.
   This sums up the environment in which the  future mathematician, who was the second of six siblings, spent his childhood.  
   
   Until age fourteen, Bernhard Riemann received his teaching at home, mostly from his father, because there was no nearby school. In 1840, his parents sent him to Hanover, where he entered secondary school, living with his grandmother. The latter died two years later, and his family moved him to the Johanneum Gynmasium in L\"uneburg,\footnote{Neuenschwander, in \cite{Neuenschwander}, says that Riemann, after spending two years at Hanover's high school, was sent by his parents to the Johanneum L\"uneburg, following difficulties with his teachers in Hanover. To support this, the author  refers to a correspondence exchanged between Riemann
and his father in 1842. The same author notes that Riemann’s certificate upon leaving
the high school in Hanover included the remark that, in recent
times, he had drawn criticism for his childish behavior
(kindisches Benehmen).}  another town in the Kingdom of Hanover,
 where he lived as a boarding pupil, from 1842 to 1846. The school had a good reputation. The inscription \emph{doctrinae, virtuti, humanitati}  (To learning, to virtue, and to humanity) hung above its entrance.
 The focus of instruction was on the classical languages. Riemann's mother, Charlotte Ebell, died in 1846.

 Initially, the young Bernhard was meant to become, like his father, a Protestant pastor. Encouraging him in this direction was a natural choice for his parents, as it would have been a way to ensure that their son would have a livelihood and a social standing. Those who knew him at that time, including his high school teachers, said that he  developed a special interest in mathematics at an early age, though it was not his unique favorite subject; he enjoyed studying in general and was gifted in all areas, especially the classics.

Since his childhood, Bernhard Riemann was a lonely boy, subject to depression and with difficulties to interact with his schoolmates. He used to return home whenever he could, to be with his family; it was the only place where he felt comfortable. But those opportunities were not very frequent, since traveling was expensive and complicated.

Riemann liked history. In a letter to his father, dated May 3, 1840
(he was 14), he complains about the fact that at the Gymnasium there were fewer
lessons on history than on Rechnen (computing).
On August 5, 1841, he
writes, again to his father, that he has the highest grade, among his classmates, in history. Besides history, Riemann was doing very well in Greek, Latin and German composition
(letters of February 1, 1845 and March 8, 1845). According to another letter to his
father, dated April 30, 1845, it is only at that point that he started being seriously 
attracted by mathematics. He mentions that his mathematics teacher, Constantin Schmalfuss, seeing his capabilities, advised him to pursue studies in mathematics.

 Schmalfuss was also the Director of the L\"uneburg Gymnasium. When he noticed how easily his pupil was making progress in mathematics,  and the degree of maturity he had attained in a short time, he decided to
exempt him from taking his class and allowed him to borrow advanced mathematics books from his own library and to read them on his own.  Riemann took this opportunity and started studying these books. This is how he became familiar at an early age with some foundational mathematical texts, including works of great mathematicians of the past
such as Euclid, Archimedes, Descartes, and Newton. On one occasion he borrowed a nine hundred page book by Legendre on number theory which, it is said, he read in six days.   
  
In the same letter to his father, April 30, 1845, Riemann declares that he plans to enroll
the University of Göttingen to study theology, but that in reality he must think more, before he decides what to do, since, if he does not take the right way, he ``will bring nothing good to any subject." The last sentence, the will of ``bringing something good to a subject", is remarkable, from a high-school pupil, and it shows how mature he was at age 19.

On April 25,1846, Riemann enrolled at the faculty of theology of the University of Göttingen, with the encouragements of his father. While he was there, he managed to attend lectures on mathematics. In particular, during the winter semester 1846/47, he followed a course by Carl Friedrich Gauss and another one by 
Moritz Stern\footnote{Moritz Abraham Stern (1807-1894), who was teaching at G\"ottingen, had earned his doctorate at that same university, before a jury chaired by Gauss. A class of prime numbers is named after him.} on definite integrals. For each course, he had to pay some amount of money. Universities in Europe were not tuition-free, as they are today.
In a letter to his father dated November 15, 1846, Riemann writes \cite{Riemann-Letters}: ``I have enrolled in six courses, for a total cost of three Louis d'or and three Thaler.\footnote{This is a mixed gold and silver coinage used in past accounting. If the transcription in \cite{Riemann-Letters} is correct, this means that the tuition  at the University of Göttingen was very costy.} Namely, Cultural History of Greece and Rome, by Hermann, Introduction to the Study of Theology by Ehrenfeuchter, Recent Church History by Reincke, General Physiology by Lotze, The Definite Integral by Stern, and the Theory of least squares by Gauss."  He adds: ``General natural history would be very interesting, and I would certainly attend, if along with everything else I had enough money."

Neuenschwander, in \cite{Neuenschwander}, reports that
the University of Göttingen has kept a list of books that
Riemann borrowed from its Library. These records show that Riemann, at that time, was interested not only in
Gauss’s works but also, and above all, in works of some
French mathematicians. He borrowed Cauchy’s \emph{Course in Analysis},
\emph{Differential Calculus}, and \emph{Mathematical Exercises},
and Legendre’s \emph{Treatise on Elliptic Functions and Eulerian integrals}.\footnote{This treatise by Legendre was published in three volumes between 1825 and 1828. The author establishes there a classification and conducts an  in-depth study of elliptic integrals and functions, laying important groundwork for this topic and for the subsequent works of  Abel, Jacobi, and Riemann on it.}  

After this, Riemann asked his father the permission to transfer to the faculty of philosophy so that he could study mathematics, and his father responded positively. And thus began the brief but dazzling career of the man who arguably is the most brilliant mathematician of modern times.

  Between 1847 and 1849, with a
break in 1848,\footnote{Major disturbances took place in Berlin in 1848, linked to a popular uprising aimed at demanding a constitution, freedom of the press, and German unity. The army fired on the crowd, triggering intense street fighting and the erection of barricades. The disturbances continued into 1849. We learn from \cite{Neuenschwander} that in the spring of 1849, yielding to his father’s will, who
was alarmed by the unrest in Berlin, Riemann returned to Göttingen,
where the situation was calmer.
}  Riemann studied mathematics at the University of Berlin. His professors included Jakob Steiner, Charles Gustave Jacob Jacobi, Johann Peter Gustav Lejeune Dirichlet and Gotthold Eisenstein, that is, some of the greatest German mathematicians of that generation. 

It appears, from the thoughts Riemann noted down, that Dirichlet was  the teacher he admired the most. Johann Peter Gustav Lejeune
 Dirichlet (1805-1859), who taught in Berlin from 1828 to 1855, is considered   the founder of analytic number theory, with  fundamental contributions to trigonometric series (not unrelated to number theory),  and to potential theory. These are all topics that were dear to Riemann. 
 Later, Riemann met Dirichlet again in Göttingen, where Riemann became a student and Dirichlet served as a professor.
After Dirichlet passed away in 1859, Riemann succeeded him as a professor and head of the mathematics department at the University of Göttingen.
 
%
%
%
%
%
  Dirichlet’s lectures in Berlin, at that epoch, were centered mostly on theoretical
physics, with a stress on partial differential equations. It is from these lectures that Riemann became
familiar with potential theory, a topic which was to play an important role in
  his later work, especially in his work on uniformization. Klein, in his \emph{Development of mathematics in the 19th century} \cite[p. 234-235]{Klein-Development},  published posthumously in 1926 and 1927, writes the following:
  \begin{quote}\small
Dirichlet loved to make things clear to himself in an intuitive substrate; along with this he
would give acute, logical analyses of foundational questions and would avoid long computations
as much as possible. His manner suited Riemann, who adopted it and worked according
to Dirichlet’s methods.
  \end{quote}
   
  Riemann’s admiration for Dirichlet is expressed at several places in his writings,
for instance in the third section of the historical part of his Habilitation dissertation
on trigonometric series which we already mentioned. Speaking of trigonometric series, Riemann describes Dirichlet as the mathematician who laid a solid foundation for a great deal of important analytical research. He recalls that Dirichlet highlighted an error made by Euler and succeeded in clarifying a question that had occupied so many eminent geometers over the preceding seventy years.  He clearly explains how his own work draws inspiration from Dirichlet’s and builds upon it. He then explains how Fourier series, in addition to their applications in physics, are also used in pure mathematics --- more specifically, in number theory --- and that Dirichlet was a pioneer in this field as well.

During his stay at Berlin, Riemann also attended lectures on optics by the physicist Wilhelm Dove (1803-1879),
which he found very interesting, and on astronomy, by the prominent German astronomer Johan Franz Enke (1791-1865)  (letter without date,
quoted in \cite{Riemann-Letters}). About the latter, Riemann says that “his presentation for the most
part is rather dry and boring, however the time that we spend at the observatory once
a week, from 6pm to 8pm, is useful and instructive.”

Besides Dirichlet, one may mention Jacobi. 
   In a letter to his father,
dated May 30, 1849, Riemann writes  \cite{Riemann-Letters}: ``Jacobi has just
begun a series of lectures in which he leads off once again with the entire system
of the theory of elliptic functions in the most advanced, but elementary way." In
another letter (written around the same period but without date), also sent from Berlin, Riemann writes: ``I enrolled with
five other students into a private class (\emph{Privatissumum}) with Eisenstein, who was
promoted in the course of this semester to a Privatdozent, with a paper on the theory
of elliptic functions."

In the fall of 1849, after two years spent in Berlin, Riemann returned to Göttingen, where he began working on his doctoral dissertation under Gauss's supervision. But the first contact between Gauss and Riemann had taken place before 
Gauss officially became Riemann’s mentor. Already in 1846, in a letter dated
November 5, Riemann informing his father about the courses he
plans to follow, mentions a course by Gauss on the method of least squares.\footnote{This is a method 
for minimizing the sum of the squared residuals, that is, the differences between observed values and the values predicted by the model. Gauss developed the method of least squares around 1794-1795 in his work on tracking planetary and stellar orbits. }
After that, and during his two year stay in Berlin (1847–1849), Riemann continued
studying thoroughly Gauss’s papers. In a letter to his father, dated May 30, 1849,
he writes:
\begin{quote}\small
Dirichlet has arranged to me to have access to the library. Without his assistance, I fear there
would have been obstacles. I am usually in the reading room by nine in the morning, to read
two papers by Gauss that are not available anywhere else. I have looked fruitlessly for a long
time in the catalog of the Royal library for another work of Gauss, which won the Copenhagen
prize,  and finally just got it through Dr. Dale of the Observatory. I am still studying it.\footnote{Gauss's paper that won the Copenhagen
prize, is titled \emph{Allgemeine Auflösung der Aufgabe: die Abbildung des Einen auf eine andere Fläche so einzurichten, dass die Abbildung dem Abgebildeten in den kleinsten Theilen ähnlich wird} (General solution of the problem: to map one surface onto another so that the image is similar to the original in its smallest parts)   \cite{Gauss-Copenhagen}. Itwas published in 1825 in the  \emph{Astronomische Abhandlungen}  (Memoirs on astronomy). This paper won a prize for a competition proposed by the Copenhagen Royal Society of Sciences in
1822 whose subject was: ``To represent the parts of a given surface onto another
surface in such a way that the representation is similar to the original in its infinitesimal parts." In this paper, Gauss proved his famous result saying (in modern terms) that on a real-analytic surface, every sufficiently small neighborhood of an arbitrary point can be mapped conformally onto a subset of the plane. Gauss did not solve the general problem of mapping conformally an arbitrary finite portion of the surface. This question was settled by Riemann, after Gauss, in his doctoral dissertation. The so-called Riemann mapping Theorem asserts that 
any non-empty, simply connected open set in the complex plane which is not the entire plane can be mapped  one-to-one conformally onto the open unit disk.}

\end{quote}

Regarding his written production, Riemann endorsed Gauss’ principle: \emph{pauca sed
matura} (few but ripe).

During the same stay in Berlin, Riemann followed lectures by Dirichlet on topics
related to Gauss’s works. He writes to his father (letter without date): ``My own course of specialization is the one with Dirichlet; he lectures on an area of mathematics
to which Gauss owes his entire reputation. I have applied myself very seriously to
this subject, not without success, I hope."

In Göttingen, after his return in 1849, Riemann attended lectures and seminars by the newly hired physics professor
Wilhelm Eduard Weber (1804-1891),  who was also Gauss’s collaborator and close
friend. Klein writes, in the \emph{Development of mathematics in the 19th century} \cite[p. 235]{Klein-Development}: 
\begin{quote}\small In Weber, Riemann found a patron and a fatherly friend. Weber recognized
Riemann’s genius and drew the shy student to him. [...] Riemann’s interest in the
mathematical treatment of nature was awakened by Weber, and Riemann was strongly
influenced by Weber’s questions.
\end{quote}

 Riemann's doctoral dissertation was defended at the University of  G\"ottingen  in 1851. This dissertation, titled \emph{Grundlagen f\"ur eine allgemeine Theorie der Functionen einer veränderlichen complexen Grösse} (Foundations for a General Theory of Functions of a Complex Variable) \cite{Riemann-Doctoral}, is probably the most important piece of work to appear as a doctoral thesis. It was submitted on December 16, 1851.  In a letter to his brother, dated
November 26, 1851, Riemann writes that after he delivered his doctoral dissertation manuscript,  Gauss took it home and examined it for a few days, but that before reading
it, Gauss told him that  for years he had been preparing an essay, on which he is still
occupied, whose subject is the same, at least in part, as that covered by this dissertation. He also told him that he had already expressed his intention to
write on that subject
  in his doctoral dissertation which, at that time, was already 52 years old.

After Riemann obtained his doctoral degree, and upon Gauss's recommendation, he was admitted, again at Göttingen University, to work on his Habilitation, the degree which would allow him to be a Privatdozent, that is, to teach and examine students at university.\footnote{Originally, the word \emph{Privatdozent}  refers to an unpaid lecturer who received fees from the students attending his
lectures.} He spent thirty months working on his Habilitation dissertation. This is the memoir which we already mentioned, called \emph{\"Uber die Darstellbarkeit einer Function durch eine trigonometrische Reihe} (On the representability of a function by a trigonometric series) \cite{Riemann-Habilitation-memoir}. In this memoir, Riemann
obtained the conditions that we now call the Riemann integrability condition, that ensure that a function has an integral. 

For his Habilitation, which he completed under Gauss, Riemann  was required,  as part of his qualifying examination, to deliver a Habilitation lecture. It was customary at G\"ottingen that the candidate prepares three topics for that lecture. Gauss (being Riemann's mentor) had to choose one of them to hear him lecturing. Riemann prepared two lectures on mathematical physics (mainly electricity and magnetism) and one on geometry. Gauss chose the lecture on geometry, and Riemann delivered his lecture, \emph{\"Uber die Hypothesen welche der Geometrie zu Grunde liegen} (On the Hypotheses Which Lie at the Foundations of Geometry), on June 10, 1854. This lecture is considered to be the founding stone of  Riemannian geometry and, beyond that, of all modern geometry, and it revolutionized mathematics.
 It consisted of two parts. In the first one, Riemann introduced the notion of $n$-dimensional space, ultimately establishing the foundation for what is known today as a Riemannian manifold.
 
 It is interesting to note that  the German  word \emph{Mannigfaltigkeit}  which Riemann chose to denote a manifold in the mathematical sense, was already used in the German philosophical literature (but not in mathematics),   in particular by Kant. Indeed, in the latter's \emph{Kritik der reinen Vernunft} (Critique of Pure Reason), 
 the word Mannigfaltigkeit, usually translated into English by
``manifold" or ``diversity", 
  designates the raw, unorganised and fragmented data (colors, sounds, textures, and shapes) that our senses receive from the external world before our minds organize it into coherent experiences. Leibniz, who usually wrote in French or Latin (and not in German, which was his mother tongue), in his famous work, \emph{Monadologie} (Monadology),   originally written  in French,
 used the word \emph{multiplicité} (multiplicity), to describe a complex plurality co-existing within a single, simple substance, a word which 
his German  translators and commentators transformed again into  \emph{Mannigfaltigkeit}. The German philosopher Johann Friedrich Herbart whom Riemann --- according to his own account, as we shall recall below --- considered himself the heir, in his philosophical writings (late 18th and early 19th centuries), argued that our minds actively construct the concept of space through psychological experience, and used extensively in this context the same word, \emph{Mannigfaltigkeit}, to mean a continuous, structured system of mental concepts.  We refer the interested reader to the article \emph{Philosophy behind mathematics} \cite{Ohshika} by K. Ohshika on the use of the term
\emph{Mannigfaltigkeit} by Riemann, and how the modern precise notion of manifold was later made out of it, through the works of Poincaré, Hilbert, Weyl and others.  On the same subject, we also refer to the articles \emph{Mathematical practice as philosophy,
with Galois, Riemann, and Grothendieck} by Arkady Plotnitsky \cite{Plotnitsky-Comprehending} and \emph{Le concept de variété et la nouvelle géométrie de l’espace dans la pensée de Bernhard Riemann. L’émergence d’une nouvelle vision des mathématiques et de ses rapports avec les sciences fondamentales} (The concept of manifold and the new geometry of space in Bernhard Riemann’s thinking. The emergence of a new vision of mathematics and its relations with fundamental sciences) \cite{Boi} by Luciano Boi.

It is also interesting to recall that Riemann was familiar, from his religious background and theological studies, with the word \emph{Dreifaltigkeit}, the  German word for Trinity,  etymologically suggesting a kind of ``three-foldedness".   The word Mannigfaltigkeit is obtained from the word Dreifaltigkeit by replacing the prefix Drei (meaning three) by the prefix Mannig (which means many).
 
    On  February  5, 1854, 
right after the submission of his  Habilitationsschrift, Riemann wrote  to his father: 
\begin{quote}\small
 I resumed my investigations into the coherence of the laws of Nature and got so involved in it that I could not tear myself loose. The continuous preoccupation with it has become bad for my health. In fact, right after New Year's, my usual affliction set in with such persistence that I could only obtain relief through the strongest remedies. As a result I felt very ill, felt unable to work, and sought to again put my health in order through long walks.
\end{quote}
This letter indicates that Riemann's illness was exacerbated by the heavy workload he placed on himself.
In a letter to his brother dated June 26, 1854, he writes, in the same tone:

\begin{quote}\small
 [\ldots] I had completed my habilitation paper at the beginning of December,
submitted it to the dean, and soon after once again turned to my investigation on the coherence of the fundamental laws in physics; also that I
so immersed myself in it that when the theme for my examination lecture
was posted at the colloquium, I could not immediately tear myself away.
Right after, I came down sick, partly, of course, as a result of too much
brooding, and partly as a result of sitting a lot in my room during bad
weather.
 \end{quote}

Gauss exerted a profound influence on Riemann. Some historians of mathematics argue that because Gauss was elderly, in poor health and notoriously averse to teaching when Riemann enrolled as his doctoral student, his guidance must have been limited. This assumption, however, contradicts the sheer scope and variety of Riemann’s mathematical ideas, many of which he himself explicitly attributed to Gauss’s direct mentorship or writings.  On this subject, there is, unfortunately, a vast body of literature written by historians of mathematics who are not mathematicians themselves and who fail to understand that, ultimately, a mathematician's influence on his students cannot be measured by the time spent in conversation. Articles that may mislead the reader continue to appear on this topic; see, for example, the recent article \cite{Rowe}. Gauss passed away in 1855, just a year after Riemann completed his Habilitation, but his intellectual imprint on Riemann remained permanent.

 Dirichlet moved to G\"ottingen in  1855, and Riemann was already there. Felix Klein writes in the work we already mentioned \cite{Klein-Development}:
\begin{quote}\small
Riemann was bound to Dirichlet by the strong inner sympathy of a like mode of thought. Dirichlet loved to make things clear to himself in an intuitive substrate; along with this he would give acute, logical analyses of foundational questions and would avoid long computations as much as possible. His manner suited Riemann, who adopted it and worked according to Dirichlet's methods.
\end{quote}
 
Riemann became (what we would call today) an associate
professor at G\"ottingen in 1857 and a full professor in 1859.

During the year 1860, Riemann visited Paris for a few weeks, and he met there, among others, the French mathematicians
Joseph Bertrand, Charles Briot, Charles Hermite, Victor Puiseux, and Joseph Serret. Among them, Briot besides being a mathematician, was a physicist and an astronomer. In particular, he devised a formula for light dispersion that bears his name. In mathematics, he worked on Abelian and elliptic functions, a subject which Riemann cherished.
Riemann wrote to his sister that he had known and admired Briot
for quite some time already, ``because of his excellent work".

In June of 1862,
Riemann married Elise Koch, a friend of his sisters. The same year, he suffered from an attack of pleurisy. Always in poor health, he spent
a large part of the years that followed in Italy,  where his doctors had advised him to go because the climate there was favorable to his health. He found financial support  and took a leave of absence from G\"ottingen for these trips to Italy, where he had mathematical colleagues with whom he could discuss.  
One of these colleagues was Enrico Betti (1823-1892), known to all topologists for the Betti numbers associated with a topological space (the name was given by Poincaré). Betti   was also a physicist (we are familiar with the Maxwell--Betti reciprocity theorem) and a specialist in potential theory. He was also among the first to understand the importance of Galois theory.  At Betti's suggestion, Riemann was offered a position at the University of Pisa, but he declined it, precisely because of his health. But he spent two consecutive summers in Pisa.  André Weil published an article titled \emph{Riemann, Betti and the birth of topology} \cite{Weil-Betti} which contains the transcription of two 
letters from Betti to his colleague and friend Placido Tardy, dated October 6 and 10, 1863,  in which he describes in detail his fascinating 
conversations with Riemann on the birth of ``analysis situs", the word used for the field that became later known as ``topology", a field of which Riemann was one of the major founders.

As soon as he returned to Göttingen, Riemann fell ill again, as his body could not tolerate the cold weather. His only child, a  daughter called Ida, whom he named after his older sister, was born in
Italy.

 In early June 1866, Riemann again decided to go
to Italy. He wanted to spend a few peaceful months on the shores of Lago Maggiore to recover his forces and finish some  works he had started. His departure was delayed a little because of the war; indeed, these were the first days of the war between Prussia and Austria.\footnote{The conflict started essentially over a dispute for territory. The key opponents were, on the one side, Austria and its allied German states, including Bavaria, Saxony, Hanover, and W\"urttemberg, and on the other side, Prussia and allies, including the Kingdom of Italy.} Riemann reached the Lago Maggiore on June 16, 1866. That was his last trip to Italy. He died there a few days later, on July 20, in Selasca on Lago Maggiore, a hamlet, like the place he was born; he is buried in the nearby cemetery there.

Riemann's  life was short, and the same was true for most people in his close family. Along with his mother, three of his sisters and one brother, died before him, deeply shadowing his life.
He published his first work at
the age of 25 and he died at the age of 39. Thus, his productive period lasted only
15 years. His collected works, which have had a profound and everlasting influence on mathematics, stand in a single slim volume.

During his lifetime, Riemann found it difficult to have friends. The only place where he felt at home was his family's village of Quickborn, where he often used to go, until his father's death in 1855. 
%
%
%
%
%
%
%
\section{A a short summary of a long story}

It is impossible to provide a serious account of Riemann’s mathematical work within the limited space of a few pages. 
Instead, in the next section, I will focus on a single topic that is dear to me personally: Riemann surfaces and their moduli. But before that, and since this is a celebration of Riemann's work, I feel committed to make at least a short list of some of his major contributions.

The reader already knows, or has understood by now, that Riemann worked on a variety of topics and that his writings had a wide and profound  impact on mathematics. In fact, Riemann's writings had a permanent influence  on philosophy and physics as well, and beyond this, on human knowledge. Besides the theory of Riemann surfaces and their moduli, which is the subject of the next section, I mention the following fields on which Riemann had an everlasting influence:  differential geometry, topology, the theories of algebraic and complex functions, integration theory, number theory and the zeta function, minimal
surfaces, uniformization, the theories of
elliptic and Abelian integrals, the hypergeometric
series, trigonometric series, electromagnetism, heat propagation,   theta functions, 
Brownian motion, probability theory, quadratic forms, electricity, gravitation, electrodynamics, magnetism,  the investigations on the nature of space, and many others. Riemann's ideas  had a great impact on several 20th century philosophers, including the mathematician-philosophers about whom I will mention more  in the later sections: Hermann Weyl, André Weil and Alexander Grothendieck. But one could also talk about Riemann's influence on the mathematician-philosopher René Thom, discussed in several chapters of the book \cite{Papadopoulos-Thom}, and also on a whole generation of French philosophers who were not mathematicians, including Gilles Deleuze, Félix Guattari, Henri Bergson, Albert Lautman, Jules Vuillemin, and there are others. To talk about Riemann's influence on all of them needs numerous articles. 

Finally, let me point out that understanding the relationship between the concepts of discreteness and continuity
in mathematics and physics, and more generally the relationship between mathematics and physics, which are topics that were dear to Riemann and which are discussed in some detail in the collective book \cite{Riemann}. For the impact of Riemann on the philosophy of mathematics and physics, I would like to mention the articles \cite{Plotnitsky-Poenaru, Plotnitsky,  Plotnitsky-Comprehending,  Plotnitsky-continuity, Mathematical-Practice} by Arkady Plotnitsky.

\section{Riemann surfaces and the problem of moduli} \label{s:Riemann-surfaces}

Riemann's work on what became known as the theory of Riemann surfaces
was partly motivated by the desire to find a domain of definition for an algebraic
multi-valued complex-analytic function (like $z\mapsto \sqrt{z})$ on which this function becomes single-valued. (As a matter of fact, this is the origin of the word ``uniformization" as it was used by Riemann).\footnote{I learned from Yi Huang that the Chinese translation of uniformization is precisely ``to render single-valued".} The study of Riemann surfaces is also closely related to, and sometimes identified with, the theory of functions
of a complex variable, which includes, as a special subfield, that of algebraic functions. This is one of the far-reaching subjects of Riemann’s investigations, and its development
is one of the most important achievements of nineteenth century mathematics.

 The main works of Riemann on this theory is 
his doctoral dissertation,  \emph{Grundlagen f\"ur eine allgemeine Theorie der Functionen einer veränderlichen complexen Grösse} (Foundations for a General Theory of Functions of a Complex Variable) \cite{Riemann-Doctoral} (1851) and his memoir \emph{Theorie der Abel'schen Functionen} (The theory of Abelian functions) \cite{Riemann-Abelian} (1857).  On Riemann's dissertation, Ahlfors writes \cite{Ahlfors}:
\begin{quote}\small
This paper marks the birth of geometric function theory [\ldots]
Very few mathematical papers have exercised an influence on the later development of mathematics which is comparable to the stimulus received from Riemann's dissertation. It contains the germ to a major part of the modern theory of analytic functions, it initiated the systematic study of topology, it revolutionized algebraic topology, and it paved the way for Riemann's own approach to differential geometry. [\ldots] 
The central idea of Riemann's thesis is that of combining geometric thought with complex analysis.
\end{quote}

In Riemann's doctoral dissertation, we find for the first time
the concept of a Riemann surface,
the first uniformization theorem (the ``Riemann mapping theorem"), and
the introduction of potential theory in the theory of Riemann surfaces (the Dirichlet principle).  Ahlfors writes in \cite[p. 4]{Ahlfors}: ``Riemann virtually puts equality signs between two-dimensional potential theory and complex function theory."
We also find, in this thesis, several fundamental topological concepts, including the notions of connectedness, of simple-connectedness, of  degree of connectivity, of coverings and of genus of a surface. We also find there the first version of the Riemann--Roch theorem, the first equivalence relation between Riemann surfaces (in the form of birational equivalence),
  the first developments of the theory   of theta functions, of the Jacobian, of Abelian integrals and of the characterization of these integrals by their period matrices. Furthermore, Riemann stated there for the first time  that the ``number of moduli" of the set of equivalence classes of a closed surface of genus $g\geq 2$ is $3g-3$.    A central question, which remained open for several decades, was to give a precise meaning to this moduli count, and the effort to answer it gave rise to the birth of the theory of moduli   and  Teichm\"uller spaces, and to works on this topic by Weierstrass, Schwarz, Dedekind, Klein, Poincar\'e, Hilbert, Fricke, Koebe, Brouwer, Weyl, Torelli,  Siegel, Teichm\"uller,  Weil, Grothendieck, Ahlfors, Bers, Mumford, Deligne, 
 Thurston, McMullen,  Witten, Yau, 
 Kontsevich, Mirzakhani and many others.

Regarding Riemann's  memoir on Abelian functions, let me quote André Weil,  from a letter to his sister Simone Weil, dated March 26, 1940, published in his 
\emph{Collected Papers} 
\cite[Vol. 1, p. 244-255]{Weil-Oeuvres-1} and in his sister's \emph{Collected Works} \cite[t. VII, vol. 1]{Weil-S}
 He writes the following, talking about algebraic functions of one variable: 
 \begin{quote}\small
 It
is generally believed that there is nothing left to do on algebraic functions of one
variable, because Riemann, who discovered almost everything we know about these
functions (I am excepting the works of Poincaré and Klein on uniformization, and
those of Hurwitz and Severi on correspondences) did not leave for us any statement
of a big problem that concerns them. I am without doubt one of the most knowledgeable
persons on this subject; certainly because I had the good fortune (in 1923)
to learn it directly from Riemann’s writings, whose memoir is of course one of the
greatest things that a mathematician has ever written; there is not a single word
there that is not of consequence. [\ldots] I am not, of course, foolish enough to compare myself to Riemann, but to add even the slightest thing to Riemann is already, as they would say in Greek, to accomplish something --- even if, in doing so, one draws (whether explicitly or implicitly) on the work of Galois, Poincaré, or Artin.

\end{quote}

         It would be unreasonable to try to develop all these topics (even one of them) in a short article. Let me only  say a few more words on the concept of Riemann surface.  
         
         Riemann viewed such a surface as a branched cover over the sphere, but above all, as a new space on which one can do  function theory.
 From the beginning of his dissertation, we read: 

 \begin{quote} \small We restrict the variables $x,y$ to a finite domain by considering as the locus of the point $O$ no longer the plane itself but a surface $T$ spread over the plane. We admit the possibility \ldots that the locus of the point $O$ is covering the same part of the plane several times.
 \end{quote}

This notion, as an abstract surface, remained misunderstood for decades, even among prominent mathematicians, who continued to think of surfaces embedded in the three-dimensional Euclidean space with self-intersection, which is not the way Riemann conceived them; see, for example, Figure 4, taken from the famous treatise by Goursat and Appell, \emph{Théorie des fonctions algébriques et de leurs intégrales --- \'Etude
des fonctions analytiques sur une surface de Riemann} (Theory of algebraic functions and their integrals: A study of analytic functions on a Riemann surface) \cite{Appell-Goursat}
(1895).

%
%
%
%
%
%
%
Regarding the theory of functions of a complex variable,
Klein, in his article \emph{ Riemann und seine Bedeutung f\"ur die Entwicklung der modernen Mathematik} (Riemann and his significance for the development of modern mathematics) (1895) \cite{Klein-significance}, attributes the beginning of this theory to Cauchy. He writes: 
\begin{quote}\small
The founder of this theory is the great French mathematician Cauchy; but only later, in Germany, did this theory assume its modern aspect which has made it the central point of our present views of mathematics. This was the result of the simultaneous efforts of two mathematicians whom we shall have to name together repeatedly, -- of Riemann and Weierstrass.
\end{quote} 
 It is interesting to see that Klein writes after this,  in a footnote: ``In the text I refrained from mentioning Gauss, who being in advance of his time in this and in other fields, anticipated many discoveries without publishing what he had found. It is very remarkable that in the papers of Gauss we find occasional glimpses of methods in the theory of functions which are completely in line with the later methods of Riemann, as if  unconsciously a transfer of leading ideas has taken place from the older to the younger mathematician".

 In his work on Riemann surfaces, Riemann made the relation between Riemann surfaces and  potential theory. Let me quote Klein about this, from the work we just mentioned \cite{Klein-significance}:
\begin{quote}\small
It should also be observed that potential theory, which in our day, owing to its
importance in the theory of electricity and in other branches of physics, is quite universally
known and used as an indispensable instrument of research, was at that time in its infancy.
It is true that Green had written his fundamental memoir as early as 1828; but this paper
remained for a long time almost unnoticed. In 1839, Gauss followed with his researches. As
far as Germany is concerned, it is mainly due to the lectures of Dirichlet that the theory was
farther developed and became known more generally; and this is where Riemann finds his
base of operation.
\end{quote}

  Klein  talks here about the importance of potential theory in electricity and in other branches of physics. Riemann never separated physics from mathematics, and in fact, he never separated  mathematics from physics and philosophy. This is what I want to communicate in the next section.  To end the present section, let me quote Klein again, from an address he delivered in Vienna on September 27, 1894:  
 ``Riemann was one of these out-of-the way scholars who let their profound thoughts ripen, silently, in their mind."

%
%
%
%
%
%

 \section{Mathematics, physics, and philosophy}

  Riemann's work had a profound and permanent impact, not only on mathematics, but also on mechanics, electromagnetism, the theory of heat, electricity, acoustics, and other fields of physics. There is also a collection of articles by Riemann on philosophy, psychology, metaphysics, and gnoseology (the philosophical study of the origin, nature, scope, and limits of human knowledge), edited by the mathematician Heinrich Weber (1842-1913) in Riemann's \emph{Collected Works} \cite{Riemann-Gesammelte}.
In this section, I will discuss briefly some ideas from physics and philosophy, and  sometimes make relations with mathematics.

For Riemann, mathematics and physics cannot be separated. His mathematical papers contain ideas that stem from physics. Many of the problems he addressed were motivated by physics, and on several occasions he employed lines of reasoning derived from physics. His work had a definite impact on this field.  Hermann Weyl and Albert Einstein, to mention only two names that are familiar to every mathematician, based their physical theories on Riemann’s mathematical work. Klein, in his
address
delivered at the 1894 general session of the Versammlung Deutscher Naturforscher
und \"Artzte (Meeting of the German naturalists and physicians) which was held in Vienna, on September
27, 1894, cf. \cite{Klein-Riemann}, writes: 

\begin{quote}\small
I must mention, first of all, that Riemann devoted much time and thought
to physical considerations. Grown up under the tradition which is represented by the combinations of the names of Gauss and Wilhelm Weber,
influenced on the other hand by Herbart’s philosophy, he endeavored
again and again to find a general mathematical formulation for the laws
underlying all natural phenomena. [\ldots] The point to which I wish to
call your attention is that these physical views are the mainspring of
Riemann’s purely mathematical investigations.
\end{quote}

Besides physics, philosophy is always in the background of Riemann's thought.  For him, mathematics, physics, and philosophy constitute three facets of the same body of knowledge.
 His Habilitation lecture \emph{\"Uber die Hypothesen, welche der Geometrie zu Grunde
liegen} (On the hypotheses that lie at the bases of geometry)  \cite{Riemann-Habilitation-Lecture} is a merging of philosophy,
geometry and physics. The fundamental questions that he explicitly addresses there are concerned with space, form, dimension, magnitude, the infinite and the infinitesimal,
the discrete and the continuous. These are precisely the questions that obsessed the ancient Greek
philosophers, starting with the Milesians and the Pythagoreans, and passing through
Plato, Aristotle, Archimedes and the other mathematicians-philosophers from Antiquity.
  Klein, in his
\emph{Development of mathematics in the 19th century} \cite[p. 235]{Klein-Development}, reporting on some
of Riemann’s ideas from what became known as his \emph{Nachlass},\footnote{This is the collection of manuscripts, notes and
correspondence that he left after his death in 1866.} writes:
 \begin{quote}\small
Riemann thinks of space as being filled with continuous matter [Stoff], which transmits the
effect of gravity, light, and electricity. He has throughout the idea of a temporal extension
of process. A remark on this topic is found in a personal letter from Gauss to Weber --- with
an express request for complete secrecy. And now I again ask, how did these things come
to Riemann? It is just mystical influence, which cannot be defined and yet cannot be clearly
grasped, of the general atmosphere of a receptive spirit.
 \end{quote}
 
 From a purely mathematical perspective, the question of whether space is continuous or discrete has had enormous consequences: it led, for instance, to the foundation of arithmetic geometry. We will return to this topic at the end of this section, where we discuss the work of Grothendieck, as well as in the final section of this article, in connection with the work of André Weil and, once again, that of Grothendieck. This ``analogy" between, on the one hand, the continuous and the discrete, and on the other hand, analysis and arithmetic, is one of the fundamental ideas expressed by André Weil in his famous letter to his sister, dated March 26, written from prison  \cite[Item 1940a]{Weil-Analogy}, which we mention without being able to say more about it due to lack of space.

It follows from Riemann's correspondence that as a student, in Göttingen and Berlin, he attended more
courses and seminars on (theoretical and experimental) physics than on mathematics, and this interest in physics remained constant during his lifetime. Since his early
twenties, he tried to develop a theory that would unify electricity, magnetism, light
and gravitation --- the same quest that Poincaré, Lorentz and Einstein addressed after him,
culminating in the theory of general relativity. The content of one of his papers, \emph{Ein Beitrag zur Electrodynamik } (A contribution to electrodynamics)\cite{Riemann-Beitrag}, published posthumously, is part of his search for the unification of the
various forces of nature.  One may also mention in this respect his essay titled \emph{Gravitation und Licht} (Gravitation and light)  \cite[p. 532-538]{Riemann-Gesammelte}, whose subject is the connection between gravitation and light. In his lectures \emph{Schwere, Electricität und Magnetismus} (Gravity, electricity and magnetism), edited
by his student Karl Hattendorf in 1880 \cite{Riemann-Schwere}, he establishes
 necessary and sufficient conditions under which Hamilton’s principle  is satisfied. This is the problem of
 motion of a free system subject to forces that depend on
 position and motion but that are time-independent. \'Emile Picard, in Chapter VII of his famous
\emph{Traité d’analyse} (Treatise on analysis) \cite{Picard-Traite}, has a chapter called \emph{Attraction and potential}. He uses there 
  (p. 167) a transformation from Riemann’s posthumous
memoir \emph{Gravity, electricity and magnetism} which we just mentioned. Let us also note that the physicist and mathematician James Clerk Maxwell discusses
Riemann’s theory of electrodynamics in his \emph{Notiz über die elektromagnetische Theorie des Lichts} (Note on the electromagnetic
theory of light), an appendix to his paper \cite{Maxwell-Method}. In the paper \emph{Riemann on Geometry, Physics, and Philosophy: Some
Remarks} \cite{Gray}
Jeremy Gray, discusses the general question of the interaction between the
mathematical ideas expressed in Riemann’s Habilitation lecture \emph{On the hypotheses that lie at the bases of
geometry} and physics and philosophy.
  
 In a note Riemann wrote about his own work \cite[p. 233]{Klein-Development}, we read: 
 \begin{quote}\small
My main work concerns a new conception of the known works of nature, their expression
by means of other basic concepts, whereby it became possible to use the experimental data
on the reciprocal actions between heat, light, magnetism and electricity to investigate their
connections with each other. I was led to this mainly through studying the works of Newton,
Euler, and --- from another aspect --- Herbart.
\end{quote}

The name of Johann Friedrich Herbart (1776–1841), which appears at the end of this quote, is poorly known to mathematicians, but is quoted several times by Riemann in his writings. A few words on Herbart are in order.

Herbart studied philosophy in Jena
under Johann Gottlieb Fichte (1762-1814), one of the founders of the philosophical movement known as German idealism, which had its origins in the theoretical and ethical writings of Kant. Herbart studied in Göttingen where he
received his doctorate and habilitation, and he later taught there pedagogy and
philosophy, before he was offered the chair formerly held by Kant in Königsberg.
His philosophy relies on  Leibniz's  theory of monads, another mathematician-philosopher. 
Herbart wrote in 1824 a
treatise titled \emph{Psychologie als Wissenschaft, neu gegründet auf Erfahrung, Metaphysik und Mathematik} (Psychology as a science newly founded on experience, metaphysics
and mathematics). He made use of infinitesimal
calculus in his research on psychology, which he considered as a
science based at the same time on experimentation, mathematics and metaphysics. He drew a parallel between this new field and physics, in the
way Newton conceived it. Sigmund Freud was profoundly influenced by Herbart. There are several sets of notes written
by Riemann and preserved in the Riemann archives in Göttingen, which concern the
philosophy of Herbart. These notes show that Riemann was highly influenced by
the philosopher, in respect of epistemology and the philosophy of science, and
in particular for his ideas on space. Most of all, Herbart and Riemann addressed the same questions. For instance, both of them addressed the question
of whether or not our knowledge of space, time, and motion is generated by  experience. In a serious discussion of Riemann's philosophy of space, one should mention Kant. But let me quote instead a passage from another mathematician-philosopher, Alexander Grothendieck, from his  vast autobiographical, poetic, and philosophical narrative \emph{Récoltes et Semailles: Réflexions et témoignage sur un passé de mathématicien} (Harvests and Sowings: Reflections and testimony about a past as a mathematician) \cite{Grothendieck-Recoltes}. The passage is worth quoting in its entirety because it contains a certain number of great ideas. Grothendieck writes \cite{Grothendieck-Recoltes}:\footnote{The English translation is mine.}

\begin{quote}\small

  It must have been fifteen or twenty years ago now, as I was leafing through the modest volume containing Riemann’s complete works, that I
was struck by a remark he made ``in passing." In this remark, he observes that the ultimate structure of space might well be
``discrete", and that the ``continuous" representations we form of it may constitute a simplification (perhaps
excessive, in the long run, \ldots ) of a more complex reality; that for the human mind, the ``continuous" was easier to grasp than the
``discontinuous,” and that it therefore serves as an ``approximation” for us to comprehend the discontinuous. This is a remark
of surprising insight coming from a mathematician, at a time when the Euclidean model of physical space
had never yet been questioned; in a strictly logical sense, it is rather the discontinuous that has traditionally served
as a technical approach to the continuous.

The developments in mathematics over the past few decades have, in fact, revealed a much closer symbiosis between
continuous and discontinuous structures than was imagined in the first half of this century. Be that as it may,
finding a ``satisfactory" model (or, if necessary, a set of such models that ``fit together" as satisfactorily as
possible \ldots), whether ``continuous", ``discrete," or of a ``mixed" nature --- such work will surely require great
conceptual imagination and a consummate flair for grasping and uncovering mathematical structures of a new kind.
This kind of imagination or ``flair" strikes me as rare, not only among physicists (where Einstein and Schrödinger
seem to have been among the few exceptions), but even among mathematicians (and here I speak from firsthand
experience).
To sum up, I foresee that the expected renewal (if it is yet to come \ldots ) will come more likely from a mathematician at
heart, well-informed about the major problems of physics, than from a physicist. But above all, it will require a person with ``philosophical openness" to grasp the crux of the problem. This problem is by no means technical in nature, but rather a
fundamental problem of the ``philosophy of nature."

\end{quote}

Arkady Plotnitsky, in his recent paper \emph{Continuity and discontinuity in, and
between, mathematics and physics} \cite{Plotnitsky-continuity}, considers the question of discreteness and continuity in modern mathematics and in physics, including this question in the later problems addressed by quantum field theories and quantum mechanics. Among other things, he discusses there the idea advanced by Riemann, followed by Grothendieck, that the continuous may serve as
an approximation of the discontinuous.

 In the next section, and in guise of a conclusion to this paper, I would like to talk briefly about Riemann's influence on some of the greatest mathematicians of the last century, and I chose Weyl, Weil, and Grothendieck, all three being, each one in his own way, a Riemann heir.
 
 \section{Weil, Weyl, Grothendieck and Riemann}
  
 I started this tribute to Riemann by mentioning his predecessors, and showing how he quoted them in his works.
 I can think of no better way to conclude the same tribute to Riemann than to quote some  major twentieth-century mathematicians  who were familiar with his work, who understood it deeply,  who incorporated it into their own work, and showing the way they quoted him.  I am  going to take this opportunity to say a few words about each of these mathematicians. This will be also an occasion to talk about some of Riemann's results that were not mentioned in the previous sections.

 I start with Hermann Weyl (1885-1955), who is one of the most prominent theoretical physicists of the twentieth century. He was the first to realise a  synthesis between general relativity and electromagnetism. He introduced the notion of gauge invariance and several other fundamental ideas in physics. 
 
 Weyl  was 
  always engaged with Riemann's mathematical, physical and philosophical legacy, frequently referring to him as a  bridge between geometry and physics.
  His work \emph{Raum $\cdot$ Zeit $\cdot$ Materie: Vorlesungen \"uber allgemeine Relativitätstheorie} (Space - Time - Matter: Lectures on General Relativity) \cite{Weyl-space}  shows that Riemannian geometry is the most natural, and often the necessary, framework for several branches of modern physics. To emphasize the naturalness of this field, Weyl makes a comparison between the fact of
moving from the flat Euclidean geometry to the curved Riemannian and the fact of moving from the action-at-a-distance physics to field theory. 
 He  writes in the work we just mentioned that 
 \begin{quote}\small
 the approach of understanding the world from its behavior in the infinitesimally small is the driving epistemological motive of near action physics and Riemannian geometry, but it is also the driving motive in the rest of Riemann’s magnificent life’s work, which is primarily directed towards complex function theory.
 \end{quote} 
In the same work, Weyl mentions the mathematicians and physicists who built on Riemann's ideas. He writes: 
\begin{quote}\small
After Riemann had made known his discoveries, mathematicians busied themselves with working out his system of geometrical ideas formally; chief among these were Christoffel, Ricci, and Levi-Civita. Riemann [\ldots]  clearly left the real development of his ideas in the hands of some subsequent scientist whose genius as a physicist could rise to equal flights with his own as a mathematician. After a lapse of seventy years this mission has been fulfilled by Einstein.
\end{quote}
 
In Section  \ref{s:Riemann-surfaces} of  the present essay, we have discussed  Riemann's Complex function theory and his introduction of the concept of Riemann surface. It was Hermann Weyl, in his work \emph{Die Idee der Riemannschen Fläche} (the concept of Riemann surface) \cite{Weyl-idee}, who, relying on Riemann's intuitive ideas,  formulated a precise definition of a surface equipped with an atlas, that is, built from local coordinate systems, paving the way for the modern definition of a manifold.

%
%
%
%

 From Hermann Weyl, we pass to André Weil (1906-1998), the founder of modern algebraic geometry who is also closely associated in several respects with Riemann.

 Weil started reading Riemann's works while he was a student at the \'Ecole Normale Supérieure, which he entered at the very precocious age 16. He writes in his autobiographical essay, \emph{Apprenticeship of a mathematician} \cite{Weil-Souvenirs}:  
   \begin{quote}\small That same year,\footnote{That is, the year he entered \'Ecole Normale.} I began to read Riemann.
[\ldots] During my year of instruction in philosophy,\footnote{This was André Weil's final year of high school, which he completed at age fourteen.} I had also been struck by a phrase of Poincaré’s
which expresses no less an extreme position: ``The value of civilizations lies only in their sciences
and arts." With such ideas in my mind I had no choice but to dive headlong into the works of the
great mathematicians of the past, as soon as they were materially and intellectually within my
grasp. Riemann was the first; I read his inaugural dissertation and his major work on Abelian
functions. Starting out thus was a stroke of luck of which I have always been grateful. These are
not hard to read as long as one realizes that every word is loaded with meaning: there is perhaps
no other mathematician whose writing matches Riemann’s for density.
\end{quote}
 
There is no doubt about the fact that if there is one open question that André Weil would have liked to settle, then this would be the  Riemann hypothesis. This problem,  formulated in 1859 by Bernhard Riemann in his eight-page article
 \emph{Über die Anzahl der Primzahlen unter einer gegebenen Grösse} (On the Number of Primes Less Than a Given Magnitude) \cite{Riemann-Zeta}
 states that the infinite series  $ \zeta(s)=\sum_{n\geq 1} n^{-s}$ (known after him as the Riemann zeta function) has all its nontrivial zeros on
the line $\mathcal{R}e (s) = 1/2$.

Weil spent a substantial part of his mathematical life in trying to prove the Riemann hypothesis. 
Because of the lack of space, it is not possible to enter here into more details about this statement, except for saying that the study of the properties of this function goes back to Euler, who was also the first to discover the relation between this function and the distribution of prime numbers. The reader interested in a discussion of and in references to Euler's work on this function may consult \cite[\S 6]{Papadopoulos-Euler}.

 André Weil formulated and proved an analogue of the Riemann hypothesis for algebraic curves over finite fields. Again, it is not possible to enter here in the 
 details, but we can say it in a minimal amount of words. The statement is equivalent to showing that the number of rational points on a curve defined over a number field is closely related  to the size of the field. Weil announced this result in two notes published in 1940 and 1941, and he wrote the detailed proofs in papers and books he published in 1942, 1946 and 1948, see \cite{Weil-1940, Weil-1941, Weil-1946, Weil-1948}. It is generally admitted that in so doing, Weil  established the foundations of modern Algebraic Geometry.
 
In an interview of André Weil \cite[p. 10]{Weil-Interview}, conducted in October 1991 (Weil was 86 years old), in response to the question ``Mathematically speaking, do you have any hopes or regrets?", Weil replied: ``When I was young, I hoped to prove the Riemann Hypothesis. When I got a little older, I still hoped to be able to read and understand a proof of the Riemann Hypothesis. Now, I would be quite content just to learn that a proof exists."
 
A recent article by Alain Connes starts with the following \cite{Connes}: ``The Riemann Hypothesis stands as perhaps the most famous unsolved problem in mathematics, having resisted all attempts at proof since Bernhard Riemann casually remarked in his 1859
memoir that  `sehr wahrscheinlich' (very probably) all nontrivial zeros of his zeta function lie on
the critical line $\mathcal{R}e (s) = 1/2$."

 André Weil's work on the Riemann Hypothesis for curves over finite fields is an illustration of Bernhard Riemann’s original idea from a shift between the continuous geometric and analytical setting to the discrete algebraic setting.
Another illustration of the same idea  is the famous Riemann--Roch--Grothendieck theorem, which we briefly introduce below. But first, we say a few words on the classical Riemann--Roch theorem.

The origin of the Riemann--Roch Theorem can be traced back to Riemann's existence theorem, which asserts the existence of meromorphic functions on
Riemann surfaces. Providing deeper insight, the classical version of the Riemann--Roch theorem focuses on the dimension of the space containing these functions on a compact Riemann surface, specifically, with poles of a certain order at some prescribed set of
points.  This theorem is stated as a formula that gives this dimension in terms of the surface's genus and the total order at its ramified points. Thus, the theorem provides a connection between topological structures and analytical concepts.
 There are several versions and
generalizations of this theorem, many developments and  many applications, including  the Hirzebruch--Riemann--Roch theorem which applies to higher-dimensional complex spaces, and the Grothendieck--Riemann--Roch theorem on which we shall say a few words.  Atiyah--Singer's index theorem, discovered in 1963, is also
considered as a generalization of the Riemann--Roch theorem.

Speaking of Grothendieck, I cannot help but recall his involvement in a radical ecological movement, particularly through the environmental and political group \emph{Survivre et vivre}, which he co-founded, his anti-war activism (this was at the time of the U.S. intervention in Vietnam), and his resignation from his professorship at IHES after discovering that the Institute was receiving funding from the Ministry of Defense.

  Grothendieck obtained a very general version of the Riemann--Roch
theorem, formulated in the language of categories and functors, which generalizes the Riemann--Roch theorem for
algebraic varieties of arbitrary dimension and defined over arbitrary ground fields. This was one of
Grothendieck’s major discoveries. In his autobiographical work \emph{Récoltes et Semailles} which we already mentioned \cite{Grothendieck-Recoltes}, he writes: ``The year 1957 is the one when I was led to extract the theme
`Riemann--Roch' (Grothendieck’s version) which overnight consecrated me as `great
star'." Grothendieck’s version of the Riemann–Roch theorem was the starting
point of K-theory. In the section called \emph{La vision --- ou douze thèmes pour une
harmonie} (The vision --- or twelve themes for a harmony) of \emph{Récoltes et semailles},  Grothendieck considers the theme he calls the Riemann--Roch--Grothendieck Yoga as one of the twelve themes which he describes as his
``great ideas" (grandes idées). There is also a discrete Grothendieck--Riemann--Roch version, proved by Grothendieck and whose full version for general schemes remains an open research field rather than a theorem.

  I would like to conclude this article by quoting three passages from Grothendieck's \emph{Récoltes et Semailles},\footnote{The quotations from Grothendieck in this paragraph are taken from a text titled ``A Stroll Through a Work, or The Child and the Mother,” which serves as an introductory essay to the entire \emph {Récoltes et Semailles}, in which he reflects on his life and on his spiritual experience, sharing images of the landscapes of his internal world,  of his vision of mathematics and of his ideas on the nature of discovery, and examining with a critical eye the scientific community of which he was part.} the work in which he presents his vision of what a mathematician is and what type of mathematician he himself belongs to.   In doing so,  he mentions  at several places the name of Riemann, not only for his mathematical results, but also as the prototype of a mathematician, and for his way of understanding mathematics.

The first passage is extracted from
 \S 6.4 of this work,  a section titled \emph{Dream and proof}. Grothendieck describes there Riemann as the model of a young student who can understand what topology is. He writes:
 \begin{quote}\small
 When I talk to
students about the topology of a sphere, and the shapes that result from a sphere by adding handles, 
things that do not surprise young children but confuse them because they think they know what
``math" is, the first spontaneous reaction I get is: ``But that’s not math! Math, of course,
is the Pythagorean theorem, the altitudes of a triangle, and quadratic polynomials \ldots" These students
aren’t any dumber than you or me; they are reacting the same way that, from time immemorial right up to the present day,
even all the mathematicians in the world have reacted,  except for people like Pythagoras or Riemann and perhaps five or
six others.
\end{quote}

 The second excerpt is taken from Section 2.5 of the same book, a section titled \emph{The heirs and the builder}:

 \begin{quote} \small
 As for me, I feel I belong to the lineage of mathematicians
whose spontaneous calling and joy is to constantly build
new houses. Along the way, they cannot help but also invent and
gradually shape all the tools, utensils, furniture, and instruments
required --- both to build the house from the foundation to the roofline,
and to abundantly furnish the future kitchens and workshops, and
to set up the house so that it is livable and comfortable. Yet, once everything is in place, down to the very last gutter and the very last stool, it is rare for the craftsman to linger long in these places, where every stone and every rafter bears the
mark of the hand that worked and set it. His place is not in the tranquility
of ready-made worlds, however welcoming and harmonious they may be, whether they
were arranged by his own hands or by those of his predecessors. Other
tasks are already calling him to new construction sites, driven by the imperious
urge of needs that perhaps only he can clearly sense, or (even more often)
by anticipating needs that only he can foresee. His place is
out in the open air.
 He is a friend of the wind and does not fear working alone,
for months and years --- and, if necessary, for a lifetime --- unless
a welcome replacement comes to his aid. He has only two hands, like
everyone else, of course, but two hands that at every moment instinctively know what
they must do, that shrink from neither the heaviest tasks nor the most
delicate ones, and that never tire of discovering and rediscovering those
countless things that constantly call upon them to be known. Two hands
are meagre, perhaps, for the world is infinite. They will never exhaust it! And
yet, two hands are a great deal \ldots
I, who am not well-versed in history, if I were to name mathematicians in that tradition, the names that immediately come to mind are those of Galois
and Riemann (in the last century) and that of Hilbert (at the beginning of the present 
century).

\end{quote}

Finally, in Section 4.2.3 of the same book, in a paragraph titled \emph{The end of a secret}, Grothendieck uses a poetic metaphor to describe the brilliance, purity, and immense creative potential he had discerned in one of his students, who immediately reminded him of Riemann.\footnote{Grothendieck is talking here about his student Pierre Deligne}  He writes\footnote{Since there are several versions of \emph{Récoltes et Semailles} on the internet, which are sometimes different in length, let me give the address of the link I am using here:
https://webusers.imj-prg.fr/~leila.schneps/grothendieckcircle/RetSbis.pdf}: 
\begin{quote}\small
I had the privilege of witnessing the first blossoming of a child’s enthusiasm, carrying the promise of a
vast unfolding.
The blossoming I had in mind is the
fruit of innocence, eager for knowledge and always ready to delight in the beauty of things both small and great
in this inexhaustible world, or in a particular part of that world (such as the vast world of mathematics \ldots ).
He alone possesses the power of profound renewal, whether it be the renewal of the self or that of
knowledge of the things of this world. It is this power that, it seems to me, found its fullest realization in the
modest figure of Riemann.
This true fulfillment is foreign to
contempt: contempt for
others (those we feel are far beneath us \ldots), or for things too ``small" or too obvious for
us to deign to take an interest in, or for those we deem to fall short of our legitimate expectations; or even
contempt for a certain \emph{dream}, perhaps, which speaks to us insistently of the things one professes to love. It is foreign
to contempt, just as it is foreign to the concept that fuels it.
\end{quote}

\bigskip

\noindent{\bf Acknowledgements.} I am thankful to S. G. Dani, Yi Huang, Ken'ichi Ohshika, Vijay Pathak and Martina Raffo for corrections and suggestions on a first version of this paper.
Thanks are also due to S. G. Dani who asked me to write a paper on Riemann for the Bulletin of the Indian Mathematics Consortium.

 \bigskip

\noindent{\bf Authors' address:}

\noindent Athanase Papadopoulos,
Institut de Recherche Mathématique Avancée
(Université de Strasbourg et CNRS),
7 rue René Descartes,
67084 Strasbourg Cedex France.

\noindent email: papadop@math.unistra.fr

 \end{document}